\documentclass[11pt, a4paper]{amsart}
\usepackage[dvips]{epsfig}
\usepackage{amsmath}
\usepackage{amssymb}
\usepackage{amsfonts}
\usepackage{amsthm}
\usepackage{amsbsy}
\usepackage{amsgen}
\usepackage{amscd}
\usepackage{amsopn}
\usepackage{amstext}
\usepackage{amsxtra}
\usepackage{enumerate}
\usepackage{hyperref}
\usepackage{lineno}
\usepackage{marginnote}
\usepackage{verbatim}

\newtheorem{theorem}{Theorem}[section]

\newtheorem{corollary}[theorem]{Corollary}

\newtheorem{lemma}[theorem]{Lemma}

\newtheorem{proposition}[theorem]{Proposition}
\newtheorem{remark}[theorem]{Remark}

\newcommand{\vol}{\operatorname{Vol}}

\newcommand{\var}{\operatorname{Var}}
\newcommand{\cov}{\operatorname{Cov}}
\newcommand{\corr}{\operatorname{Corr}}
\newcommand{\proj}{\operatorname{Proj}}

\numberwithin{equation}{section}
\begin{document}

\title[Nodal length of random spherical harmonics]{The Asymptotic Equivalence of the Sample Trispectrum and the Nodal Length for Random Spherical Harmonics}

\date{\today}

\author{ Domenico Marinucci, Maurizia Rossi and Igor Wigman}

\address{Department of Mathematics, Universit\`a di Roma Tor Vergata (Corresponding author)} \email{marinucc@mat.uniroma2.it}

\address{Unit\'e de Recherche en Math\'ematiques, Universit\'e du Luxembourg} \email{maurizia.rossi@uni.lu}

\address{Department of Mathematics, King's College London } \email{igor.wigman@kcl.ac.uk}

\maketitle



\begin{abstract}
We study the asymptotic behaviour of the nodal length of random $2d$-spherical
harmonics $f_{\ell}$ of high degree $\ell \rightarrow\infty$, i.e. the length of their zero set
$f_{\ell}^{-1}(0)$.
It is found that the nodal lengths are asymptotically equivalent, in the $L^{2}$-sense, to
the "sample trispectrum", i.e., the integral of $H_{4}(f_{\ell}(x))$, the
fourth-order Hermite polynomial of the values of $f_{\ell}$. A particular by-product of this
is a Quantitative Central Limit Theorem (in Wasserstein distance) for the nodal length, in the high energy limit.
\end{abstract}

\begin{itemize}
\item \textbf{AMS Classification}: 60G60, 62M15, 53C65, 42C10, 33C55.

\item \textbf{Keywords and Phrases}: Nodal Length, Spherical Harmonics,
Sample Trispectrum, Berry's Cancellation, Quantitative Central Limit Theorem
\end{itemize}

\section{Introduction and Main Results}

\subsection{Background}

Let $\mathbb{S}^{2}$ be the unit $2$d sphere and
$\Delta _{\mathbb{S}^{2}}$ be the Laplace-Beltrami operator on
$\mathbb{S}^{2}$. It is well-known that the spectrum of $\Delta _{\mathbb{S}^{2}}$
consists of the numbers $\lambda _{\ell }=\ell (\ell +1)$ with $\ell\in \mathbb{Z}_{\ge 0}$,
and the eigenspace corresponding to $\lambda_{\ell}$ is the $(2\ell+1)$-dimensional linear
space of degree $\ell$ spherical harmonics. For $\ell\ge 0$ let
$\left\{ Y_{\mathbb{\ell }m}(.)\right\} _{m=-\ell ,\dots ,\ell}$ be an arbitrary $L^{2}$-orthonormal basis
of real valued spherical harmonics $$Y_{\mathbb{\ell }m}:\mathbb{S}^{2}\rightarrow
\mathbb{\mathbb{R}}$$ satisfying
\begin{equation*}
\Delta _{\mathbb{S}^{2}}Y_{\mathbb{\ell }m}+\lambda _{\ell }Y_{\mathbb{\ell
}m}=0\text{ },\hspace{1cm}Y_{\mathbb{\ell }m}:\mathbb{S}^{2}\rightarrow
\mathbb{\mathbb{R}}.
\end{equation*}%

On $\mathbb{S}^{2}$ we consider a family of Gaussian random fields (defined on a suitable probability space $(\Omega,\mathcal{F},\mathbb{P})$)
\begin{equation}
f_{\ell }(x)=\sqrt{\frac{4\pi }{2\ell +1}}\sum_{m=-\ell }^{\ell }a_{\ell
m}Y_{\ell m}(x),  \label{fell}
\end{equation}%
where the coefficients
$\left\{ a_{\ell m}\right\} _{m=-\ell ,...,\ell }$ are i.i.d. standard
Gaussian random variables (zero mean and unit variance); it is immediate
to see that the law of the process $\left\{ f_{\ell }(.)\right\} $ is
invariant with respect to the choice of a $L^{2}$-orthonormal basis $%
\{Y_{\ell m}\}$. The random fields $\{f_{\ell }(x),\;x\in \mathbb{S}^{2}\}$ are centred,
Gaussian and isotropic, satisfying
$$\Delta _{\mathbb{S}^{2}}f_{\mathbb{\ell }}+\lambda _{\ell }f_{\mathbb{\ell }}=0;$$
these are the {\em random} degree-$\ell$ spherical harmonics.
From the addition formula for spherical harmonics \cite[(3.42)]{MaPeCUP},
the covariance function of $f_{\ell}$ is given by
\begin{equation*}
\mathbb{E}[f_{\ell }(x)\cdot f_{\ell }(y)]=P_{\ell }(\cos d(x,y)),
\end{equation*}%
where $P_{\ell }$ are the Legendre polynomials,
and $d(x,y)$ is the spherical geodesic distance between $x$ and $y,$ $%
d(x,y)=\arccos (\left\langle x,y\right\rangle )$.
The random spherical harmonics naturally arise
from the spectral analysis of isotropic spherical random fields (e.g. \cite{CMW, CW, CaMar2016}), and
Quantum Chaos (e.g. \cite{MaPeCUP, {wigmanreview}}); their geometry is of significant interest.

In this paper, we shall focus on the nodal length of the random fields
$\left\{f_{\ell }(.)\right\} ,$ i.e. the length of the nodal line:
\begin{equation*}
\mathcal{L}_{\ell }:=\operatorname{len}\left\{ f_{\ell }^{-1}(0)\right\}.
\end{equation*}%
Here $\left\{ \mathcal{L}_{\ell }\right\}_{\ell\ge 0} $
is a sequence of random variables; Yau's conjecture \cite{Yau}, asserts that the nodal volume of Laplace eigenfunctions
on {\em smooth} $n$-manifolds
is commensurable with the square root of the eigenvalue. An application of Yau's conjecture, established \cite{DonnFeffa} for all {\em analytic} manifolds, on the sample functions $f_{\ell}$ implies that
one has, for some {\em absolute} constants $C\geq c>0$%
\begin{equation}
\label{eq:Yau's conj}
c\sqrt{\lambda _{\ell }}\leq \mathcal{L}_{\ell }\leq C\sqrt{\lambda _{\ell }}%
\text{ for all }\ell\ge 1.
\end{equation}
The lower bound in \eqref{eq:Yau's conj} was recently established \cite{limonov} for all smooth manifolds.

While the expected value of $\mathcal{L}_{\ell }$ was computed \cite{Neuhaisel} by a standard
application of the Kac-Rice formula to be
\begin{equation*}
\mathbb{E}\left[ \mathcal{L}_{\ell }\right] =
\left\{ \frac{\lambda _{\ell }}{2}\right\} ^{1/2}\times 2\pi;
\end{equation*}
evaluating the variance proved to be more subtle, and was shown \cite{Wig} to be asymptotic to
\begin{equation}
\label{eq:var nod length log}
\var\left\{ \mathcal{L}_{\ell }\right\} =\frac{\log \ell }{32}+O(1).
\end{equation}%
It follows that the "generic" (Gaussian) spherical eigenfunctions obey a stronger
law than \eqref{eq:Yau's conj}, with normalised nodal length $\frac{\mathcal{L}_{\ell }}{\ell}$ converging
to a positive constant.

\subsection{Main Results}

In this work, we are interested in the analysis of the fluctuations of the
nodal length around its expected value; in particular a
(quantitative) central limit theorem will be established for the (centred and standardized)
fluctuations of $\mathcal{L}_{\ell }.$ This result is a
rather straightforward corollary of a deeper result, namely the
asymptotic equivalence (in the $L^{2}(\Omega )$ sense) of the nodal length
and the sample trispectrum of $\left\{ f_{\ell }\right\} ,$ i.e., the
integral of $H_{4}(f_{\ell }(x))$, where $H_{4}$ is the fourth-order Hermite polynomial; we recall that
$$H_{4}(u)=u^{4}-6u^{2}+3.$$ More
precisely, let us define the sequence of centred random variables
\begin{align}
\mathcal{M}_{\ell }& :=-\frac{1}{4}\sqrt{\frac{\ell (\ell +1)}{2}}\frac{1}{4!%
}\int_{\mathbb{S}^{2}}H_{4}(f_{\ell }(x))dx=-\frac{1}{4}\sqrt{\frac{\ell
(\ell +1)}{2}}\frac{1}{4!}h_{\ell ;4}\text{ , \ }  \label{Melle} \\
h_{\ell ;4}& :=\int_{\mathbb{S}^{2}}H_{4}(f_{\ell }(x))dx\text{ , }\ell
=1,2,....;  \label{Helle}
\end{align}
the sequence $\left\{ h_{\ell ;4}\right\}$ (which we call the sample trispectrum of $f_{\ell}$) was studied earlier, and indeed
building upon \cite[Lemma 3.2]{MW2014}, it is immediate to establish the
following:

\begin{lemma}
\bigskip \label{VarH4} As $\ell \rightarrow \infty ,$ we have%
\begin{equation}
\label{eq:trispec var log}
\var\left\{ \mathcal{M}_{\ell }\right\} =\frac{1}{32}\log \ell +O(1).
\end{equation}
\end{lemma}

By means of Kac-Rice formula \cite{adlertaylor,azaiswschebor,Wig}, the spherical nodal length $\mathcal{L}_{\ell }$ can be formally
written as%
\begin{equation}
\mathcal{L}_{\ell }=\int_{\mathbb{S}^{2}}\left\Vert \nabla f_{\ell
}(x)\right\Vert \delta (f_{\ell }(x))dx, \label{Kac-Rice}
\end{equation}%
where $\delta (.)$ denotes the Dirac delta function and $\left\Vert
{\cdot }\right\Vert $ the standard Euclidean norm in $\mathbb{R}^{2}$; this
representation can be shown to hold almost surely in $\Omega ,$
and it is shown in the Appendix that it also holds in $L^{2}(\Omega )$.
Denote by $\widetilde{\mathcal{L}}_{\ell }$ the standardized nodal length,
i.e.,
\begin{equation}
\widetilde{\mathcal{L}}_{\ell }:=\frac{\mathcal{L}_{\ell }-\mathbb{E}%
\mathcal{L}_{\ell }}{\sqrt{\var(\mathcal{L}_{\ell })}}\text{ .}
\label{Lellehat}
\end{equation}
Note that the variance of $\mathcal{M}_{\ell }$ is asymptotic to the one of
$\mathcal{L}_{\ell }$, i.e.
\begin{equation*}
\frac{\var\left\{ \mathcal{L}_{\ell }\right\} }{\var\left\{ \mathcal{M}_{\ell
}\right\} }=1+O\left(\frac{1}{\log \ell }\right)\text{ , as }\ell \rightarrow \infty
\text{ ;}
\end{equation*}%
we shall also standardize the zero-mean sequence $\left\{ \mathcal{M}_{\ell
}\right\} ,$ writing
\begin{equation}
\widetilde{\mathcal{M}}_{\ell }:=\frac{\mathcal{M}_{\ell }}{\sqrt{\var(%
\mathcal{M}_{\ell })}}.  \label{Mellehat}
\end{equation}
The main contribution of the present manuscript is establishing the following asymptotic
representation for $\widetilde{\mathcal{L}}_{\ell }$:

\begin{theorem}
\label{maintheorem} As $\ell \rightarrow \infty ,$ we have that%
\begin{equation*}
\mathbb{E}\left[ \left\{ \widetilde{\mathcal{L}}_{\ell }-\widetilde{\mathcal{%
M}}_{\ell }\right\} ^{2}\right] =O\left(\frac{1}{\log \ell }\right),
\end{equation*}%
and thus in particular%
\begin{equation*}
\widetilde{\mathcal{L}}_{\ell }=\widetilde{\mathcal{M}}_{\ell }+O_{p}\left(\frac{1%
}{\sqrt{\log \ell }}\right).
\end{equation*}
\end{theorem}

In other words, after centring and normalization the spherical nodal lengths \eqref{Lellehat} and the sample trispectrum \eqref{Mellehat} are asymptotically equivalent in $L^{2}(\Omega )$ (and thus in
probability and in law).
Now recall that the Wasserstein distance between two random variables $X$
and $Y$ is given by (see e.g. \cite[Appendix C]{noupebook})
\begin{equation*}
d_{W}(X,Y)=\sup_{h:\left\Vert h\right\Vert _{Lip}\leq 1}\left\vert
Eh(X)-E(h(Y)\right\vert;
\end{equation*}%
convergence in mean-square implies convergence in Wasserstein distance, and
both imply convergence in distribution. Let $\mathcal{N}(0,1)$
denote a standard Gaussian random variable; in view of the aforementioned CLT
\cite{MW2014} on $\left\{ \widetilde{\mathcal{M}}_{\ell }\right\} $, it follows directly from Theorem \ref{maintheorem}
that:

\begin{corollary}
\label{CLT} As $\ell \rightarrow \infty ,$ we have that
\begin{equation*}
d_{W}(\widetilde{\mathcal{L}}_{\ell },\mathcal{N}(0,1))=O\left(\frac{1}{\sqrt{%
\log \ell }}\right).
\end{equation*}
\end{corollary}

Hence we obtain here a new Quantitative Central Limit Theorem (in Wasserstein distance) for the spherical nodal length.

\subsection{Discussion and Overview of Some Related Literature}

Theorem \ref{maintheorem} is closely related to the recent characterization
of the asymptotic distribution for the nodal length of arithmetic random
waves, i.e. Gaussian eigenfunctions on the two-dimensional torus $\mathbb{T}%
^{2},$ which was established in \cite{GAFA2016}. The approach in the latter
paper can be summarized as follows: the nodal length can be decomposed into
so-called Wiener-Chaos components, i.e., it can be projected on the
orthogonal subspaces of $L^{2}(\Omega )$ spanned by linear combinations of
multiple Hermite polynomials of degree $q;$ more precisely, we have the
orthogonal decomposition%
\begin{equation*}
L^{2}(\Omega )=\bigoplus\limits_{q=0}^{\infty }\mathcal{H}_{q},
\end{equation*}%
where $\mathcal{H}_{q}$ denotes the $q$-th Wiener chaos, i.e., the linear
span of Hermite polynomials of order $q$, and hence
\begin{equation}
\mathcal{L}_{n}=\sum_{q=0}^{\infty }\text{Proj}[\mathcal{L}_{n}|q]
\label{Hermexp}
\end{equation}%
where $\mathcal{L}_{n}$ denotes the nodal length of Gaussian arithmetic
random waves of degree $n=a^{2}+b^{2},$ where $n,a,b\in \mathbb{N}$, and Proj%
$[.|q]$ projection on $\mathcal{H}_{q}$, see \cite{KKW,GAFA2016} for
details. For arithmetic random waves, all the terms
$\left\{ Proj[\mathcal{L}_{n}|q],\text{ }q\text{ odd}\right\} $ in
the expansion \eqref{Hermexp} vanish for symmetry reasons, and so does the term
corresponding to $q=2.$ The latter phenomenon is one interpretation of the
so-called {\em Berry's cancellation}, i.e. the fact that the nodal
length variance is of order of magnitude smaller than the natural scaling.
Indeed it has been shown \cite{GAFA2016,ROSSI2015} that
the term corresponding to $q=2$ dominates the fluctuations of the boundary
length of excursion sets for arithmetic random waves for any threshold value
$z\neq 0;$ for $z=0$ it vanishes, and the dominating term
is the projection onto the $4$th order chaos.

The asymptotic domination of the second-order chaos for $z\neq 0,$ and its
disappearance for $z=0,$ have been shown recently to occur for other
geometric functionals of excursion sets of random eigenfunctions in a
variety of circumstances, such as the excursion area and the Defect (\cite{DI,MW2014,MR2015}
covering all dimensions $d\geq 2$), and the Euler-Poincar\`{e} characteristic \cite{CaMar2016} (see also \cite{DNPR2016}).
The fact that a single chaos dominates clearly allows for a much neater
derivation of asymptotic distribution results; in particular, quantitative
central limit theorems have been given \cite{DI,MW2014,MR2015,CaMar2016}
for various geometric functionals of random spherical
harmonics, in the high-energy limit where $\lambda _{\ell } \rightarrow \infty$; on the torus
the asymptotic behaviour is more complicated, depending on the different subsequences as $n$ grows
\cite{KKW,GAFA2016,PR2017} for the nodal length of arithmetic random waves, and \cite{RW17} for nodal
intersections of arithmetic random waves against a fixed curve.

The results we shall give here confirm the asymptotic dominance of the
fourth-order component; in this sense, they are analogous to those in \cite%
{GAFA2016} for the case of the torus. On the other hand, here we are able to
obtain a neater expression for the leading term, which is of
independent interest, and makes the derivation of a Quantitative Central
Limit Theorem much more elegant. In fact, rather than studying the asymptotic
behaviour of the fourth-order chaos (which is a sum of six terms involving
the eigenfunctions and their gradients), we establish the asymptotically
full correlation of the nodal length with a term which can be evaluated in
terms of the eigenfunctions themselves, and not their gradient components.
The resulting approximation (valid in the mean square sense) is therefore
surprisingly simple, and the quantitative central limit theorem follows as an
immediate consequence. We believe that this technique can be applicable to
other examples of geometric functionals for random spherical harmonics.

\vspace{2mm}

As mentioned earlier, the approach we use in this paper does not require to
study directly the asymptotic behaviour of the full components in the
fourth-order chaos, as it was done earlier in \cite{GAFA2016} for
eigenfunctions on the torus. On the other hand it is obvious that the
random variables $\mathcal{M}_{\ell } \in \mathcal{H}_{4},$ and a
quick inspection to the proof of our main result reveals that we have also
the following asymptotic equivalence: As $\ell \rightarrow \infty ,$ we have
that%
\begin{equation*}
\corr\left\{ \proj[\mathcal{L}_{\ell }|4],\mathcal{M}_{\ell }\right\} =1+O\left(%
\frac{1}{\log \ell }\right),
\end{equation*}%
and hence%
\begin{equation*}
\mathbb{E}\left[ \left\{ \proj[\widetilde{\mathcal{L}}_{\ell }|4]-\widetilde{%
\mathcal{M}}_{\ell }\right\} ^{2}\right] =O\left(\frac{1}{\log \ell }%
\right)\Longrightarrow \proj[\widetilde{\mathcal{L}}_{\ell }|4]=\widetilde{\mathcal{%
M}}_{\ell }+O_{p}(\frac{1}{\sqrt{\log \ell }}).
\end{equation*}%
Likewise, as $\ell \rightarrow \infty ,$ we also have that%
\begin{equation*}
\var\left\{ \proj[\mathcal{L}_{\ell }|4]\right\} =\frac{\log \ell }{32}%
+O(1)\Rightarrow \frac{\var\left\{ \proj[\mathcal{L}_{\ell }|4]\right\} }{%
\var\left\{ \mathcal{L}_{\ell }\right\} }=1+O\left(\frac{1}{\log \ell }\right),
\end{equation*}%
and hence%
\begin{equation*}
\mathbb{E}\left[ \left\{ \widetilde{\mathcal{L}}_{\ell }-\proj[\widetilde{%
\mathcal{L}}_{\ell }|4]\right\} ^{2}\right] =O\left(\frac{1}{\log \ell }%
\right)\Longrightarrow \widetilde{\mathcal{L}}_{\ell }=\proj[\widetilde{\mathcal{L}}%
_{\ell }|4]+O_{p}\left(\frac{1}{\sqrt{\log \ell }}\right).
\end{equation*}%
In other words, it does follow from our results that the fourth-order chaos
projection dominates the high-frequency behaviour of the spherical nodal
length, as for the two-dimensional toroidal eigenfunctions.

As discussed before, the nodal length of random spherical harmonics can be
viewed as the special case (for $z=0$) of the boundary length of excursion
sets ($\mathcal{L}_{\ell }(z),$ say, with $\mathcal{L}_{\ell }:=\mathcal{L}%
_{\ell }(0))$. For $z\neq 0,$ it was shown in \cite{ROSSI2015}, Proposition
7.3.1 (see also \cite{CaMar2016}, Subsection 1.2.2, and \cite{GAFA2016},
Remark 2.4) that the dominant term corresponds to the second order chaos,
which can be expressed as%
\begin{equation*}
Proj[\mathcal{L}_{\ell }(z)|2]=2\left\{ \frac{\lambda _{\ell }}{2}\right\}
^{1/2}\sqrt{\frac{\pi }{8}}\left\{ z^{2}\phi (z)\right\} \frac{1}{2!}\int_{%
\mathbb{S}^{2}}H_{2}(f_{\ell }(x))dx,
\end{equation*}%
a component that vanishes identically for $z=0.$ On the other hand, for the
nodal case from the results in this paper one obtains the related expression%
\begin{equation*}
Proj[\mathcal{L}_{\ell }(0)|4]=-\left\{ \frac{\lambda _{\ell }}{2}\right\}
^{1/2}\sqrt{\frac{\pi }{8}}\phi (0)\frac{1}{4!}\left\{ \int_{\mathbb{S}%
^{2}}H_{4}(f_{\ell }(x))dx+O_{p}(\frac{1}{\ell })\right\}.
\end{equation*}%
It is instructive to compare these expressions with the results provided by
the Gaussian Kinematic Formula (see e.g., \cite{taylor}, \cite{adlertaylor})
for the expected value of the boundary length, which in terms of
Wiener-chaos projections can be written in this framework as
\begin{equation*}
\mathbb{E}\left[ \mathcal{L}_{\ell }(z)\right] =Proj[\mathcal{L}_{\ell
}(z)|0]=2\left\{ \frac{\lambda _{\ell }}{2}\right\} ^{1/2}\sqrt{\frac{\pi }{8%
}}\phi (z)\int_{\mathbb{S}^{2}}H_{0}(f_{\ell }(x))dx.
\end{equation*}%
We leave as an issue for further research to determine whether similarly
neat expressions can be shown to hold in greater generality, i.e. in
dimension greater than two, for higher-order chaos projections, or for
different geometric functionals.

\subsection{Outline of the paper}

In section \ref{L2Discussion} we
discuss some issues concerning the $L^{2}$ expansion of the spherical nodal
length into Wiener chaos components, and we present the analytic expression
for the fourth-order chaos (which corresponds to the leading
non-deterministic term); in section \ref{ProofTheorems} we give the proofs
of the two main Theorems, which are largely based on a Key Proposition whose
proof is collected in section \ref{ProofProposition}. The Appendix (section %
\ref{Appendix}) collects the justification for the $L^{2}$ expansion of the
nodal length into Hermite polynomials and some elementary facts about the
covariances of random spherical harmonics and their derivatives.

\subsection{Acknowledgements}

The research leading to these results has received funding from the European Research Council under the European Union's Seventh
Framework Programme (FP7/2007-2013) / ERC grant agreements n$^{\text{o}}$ 277742
\emph{Pascal} (Domenico Marinucci and Maurizia Rossi) and n$^{\text{o}}$ 335141 \emph{Nodal} (Igor Wigman), and
by the grant F1R-MTH-PUL-15STAR (STARS) at Luxembourg University (Maurizia Rossi).

\section{\label{L2Discussion} The $L^{2}$ Expansion of Nodal Length}

In this section, we present the Wiener chaos expansion of the nodal length $%
\mathcal{L}_{\ell }$ as in \eqref{Kac-Rice}.
The details of this derivation are similar to those
given in \cite{GAFA2016} (see also \cite{KratzLeon}). Let us first recall
the expression for the projection coefficients in the Hermite expansions of
the two-dimensional norm and the Dirac $\delta $-function. For independent,
standard Gaussian variables $\zeta ,\eta $ the expansion of the Euclidean
norms has been been established to be (see i.e., \cite{KratzLeon},
\cite{GAFA2016})
\begin{equation*}
\left\Vert (\zeta ,\eta )\right\Vert =\sum_{q=0}^{\infty }\sum_{n,m:2n+2m=q}%
\frac{\alpha _{2n,2m}}{(2n)!(2m)!}H_{2n}(\zeta )H_{2m}(\eta )
\end{equation*}%
where
\begin{equation*}
\alpha _{2n,2m}:=\sqrt{\frac{\pi }{2}}\frac{(2n)!(2m)!}{n!m!}\frac{1}{2^{n+m}%
}p_{n+m}\left(\frac{1}{4}\right),
\end{equation*}%
and $p_{N}$ is the swinging factorial coefficient
\begin{equation*}
p_{N}(x):=\sum_{j=0}^{N}(-1)^{j}(-1)^{N}\left(
\begin{array}{c}
N \\
j%
\end{array}%
\right) \frac{(2j+1)!}{(j!)^{2}}x^{j}.
\end{equation*}%
For the first few terms we have
\begin{equation*}
\alpha _{00}=\sqrt{\frac{\pi }{2}}\text{ ; }\alpha _{02}=\frac{1}{2}\sqrt{%
\frac{\pi }{2}}\text{ ; }\alpha _{04}=-\frac{3}{8}\sqrt{\frac{\pi }{2}}.
\end{equation*}%

On the other hand, the first few coefficients for the
expansion into Wiener chaoses of the Dirac delta function $\delta$-function
are given by (\cite{KratzLeon}, \cite{GAFA2016}):
\begin{equation*}
\beta _{0}=\frac{1}{\sqrt{2\pi }}\text{ };\text{ }\beta _{2}=-\frac{1}{\sqrt{%
2\pi }}\text{ ; }\beta _{4}=\frac{3}{\sqrt{2\pi }}.
\end{equation*}%
The Wiener-chaos decompositions need to be evaluated on variables of unit variance;
this requires dividing the
derivatives by $\sqrt{\ell (\ell +1)/2}\sim \ell /\sqrt{2}$ (here and
everywhere else $a_{\ell }\sim b_{\ell }$ means that the ratio between the two
sequences converges to unity as $\ell \rightarrow \infty $). The
$L^{2}(\Omega )$ expansion of the nodal length \eqref{Kac-Rice} then takes the form
\begin{eqnarray*}
\mathcal{L}_{\ell }-\mathbb{E}\mathcal{L}_{\ell } &=&\sqrt{\frac{\ell (\ell
+1)}{2}}\sum_{q=2}^{\infty }\sum_{u=0}^{q}\sum_{k=0}^{u}\frac{\alpha
_{k,u-k}\beta _{q-u}}{k!(u-k)!(q-u)!}\times \\
&&\times \int_{\mathbb{S}^{2}}H_{q-u}(f_{\ell }(x))H_{k}(\frac{\partial
_{1;x}f_{\ell }(x)}{\sqrt{\ell (\ell +1)/2}})H_{u-k}(\frac{\partial
_{2;x}f_{\ell }(x)}{\sqrt{\ell (\ell +1)/2}})dx
\end{eqnarray*}%
\begin{equation*}
=\sum_{q=2}^{\infty }\int_{\mathbb{S}^{2}}\Psi _{\ell }(x;q)dx,
\end{equation*}%
where%
\begin{equation*}
\Psi _{\ell }(x;q):=\sqrt{\frac{\ell (\ell +1)}{2}}\sum_{u=0}^{q}%
\sum_{k=0}^{u}\frac{\alpha _{k,u-k}\beta _{q-u}}{k!(u-k)!(q-u)!}%
H_{q-u}(f_{\ell }(x))H_{k}(\frac{\partial _{1;x}f_{\ell }(x)}{\sqrt{\ell
(\ell +1)/2}})H_{u-k}(\frac{\partial _{2;x}f_{\ell }(x)}{\sqrt{\ell (\ell
+1)/2}});
\end{equation*}%
here, we are using spherical coordinates (colatitude $\theta ,$ longitude $%
\varphi $) and for $x=(\theta _{x},\varphi _{x})$ we are using the notation
\begin{equation*}
\partial _{1;x}=\left. \frac{\partial }{\partial \theta }\right\vert
_{\theta =\theta _{x}}\text{ , }\partial _{2;x}=\left. \frac{1}{\sin \theta }%
\frac{\partial }{\partial \varphi }\right\vert _{\theta =\theta _{x},\varphi
=\varphi _{x}}.
\end{equation*}%
In particular, the projection of the nodal length on the fourth-order chaos
has the expression%
\begin{equation*}
\proj[\widetilde{\mathcal{L}}_{\ell }|4]=\int_{\mathbb{S}^{2}}\Psi _{\ell
}(x;4)dx
\end{equation*}%
\begin{equation*}
=\sqrt{\frac{\ell (\ell +1)}{2}}\left\{ \frac{\alpha _{00}\beta _{4}}{4!}%
\int_{\mathbb{S}^{2}}H_{4}(f_{\ell }(x))dx+\frac{\alpha _{20}\beta _{2}}{2!2!%
}\int_{\mathbb{S}^{2}}H_{2}(f_{\ell }(x))H_{2}(\frac{\partial _{1;x}f_{\ell
}(x)}{\sqrt{\ell (\ell +1)/2}})dx\right.
\end{equation*}%
\begin{equation*}
+\frac{\alpha _{40}\beta _{0}}{4!}\int_{\mathbb{S}^{2}}H_{4}(\frac{\partial
_{1;x}f_{\ell }(x)}{\sqrt{\ell (\ell +1)/2}})dx+\frac{\alpha _{22}\beta _{0}%
}{2!2!}\int_{\mathbb{S}^{2}}H_{2}(\frac{\partial _{1;x}f_{\ell }(x)}{\sqrt{%
\ell (\ell +1)/2}})H_{2}(\frac{\partial _{2;x}f_{\ell }(x)}{\sqrt{\ell (\ell
+1)/2}})dx
\end{equation*}%
\begin{equation}
\left. +\frac{\alpha _{02}\beta _{2}}{2!2!}\int_{\mathbb{S}%
^{2}}H_{2}(f_{\ell }(x))H_{2}(\frac{\partial _{2;x}f_{\ell }(x)}{\sqrt{\ell
(\ell +1)/2}})dx+\frac{\alpha _{04}\beta _{0}}{4!}\int_{\mathbb{S}^{2}}H_{4}(%
\frac{\partial _{2;x}f_{\ell }(x)}{\sqrt{\ell (\ell +1)/2}})dx\right\}.  \label{fourthchaos}
\end{equation}

\section{\label{ProofTheorems} Proof of the Main Results (Theorem \protect
\ref{maintheorem} and Corollary \protect\ref{CLT})}

\subsection{Proof of Lemma \protect\ref{VarH4}}

\begin{proof}
It is sufficient to note that%
\begin{eqnarray*}
\var\left\{ \mathcal{M}_{\ell }\right\} &=&\frac{\ell (\ell +1)}{2\times
4^{2}\times 24^{2}}\times \left\{ 8\pi ^{2}\int_{0}^{\pi }\mathbb{E}\left[
H_{4}(f_{\ell }(\overline{x}))\right] H_{4}(f_{\ell }(y(\theta )))\sin
\theta d\theta \right\} \\
&=&\frac{\ell (\ell +1)}{2\times 4^{2}\times 24^{2}}\times 576\frac{\log
\ell }{\ell ^{2}}+O(1)=\frac{1}{32}\log \ell +O(1),
\end{eqnarray*}%
where we have used the asymptotic result \cite[Lemma 3.2]{MW2014}
\begin{equation*}
\mathbb{E}\left[ \left\{ \int_{\mathbb{S}^{2}}H_{4}(f_{\ell }(x))dx\right\}
^{2}\right] =576\frac{\log \ell }{\ell ^{2}}+O\left(\frac{1}{\ell ^{2}}\right)\text{ ,
as }\ell \rightarrow \infty.
\end{equation*}%
\end{proof}

Before we proceed with the proof, we need to introduce some more notation. We shall write $\overline{x}=(0,0)$ for the "North Pole" and$\ y(\theta
)=(0,\theta )$ for the points on the meridian where $\varphi =0.$ Also, we transform variables as
\begin{equation*}
\psi :=L\theta \text{ , for }L:=\left(\ell +\frac{1}{2}\right), \text{ whence } y(\theta
)=y(\frac{\psi }{L});
\end{equation*}%
finally, let us define the following \emph{2-point cross-correlation
function}
\begin{equation}
\mathcal{J}_{\ell }(\psi ;4):=\left[ -\frac{1}{4}\sqrt{\frac{\ell (\ell +1)}{%
2}}\frac{1}{4!}\right] \times \frac{8\pi ^{2}}{L}\mathbb{E}\left\{ \Psi
_{\ell }(\overline{x};4)H_{4}(f_{\ell }(y(\frac{\psi }{L})))\right\}. \label{CrossCorrelation}
\end{equation}%

Our main result will follow from the following Key Proposition:

\begin{proposition}
\bigskip \label{Key} For any constant $C>0$, uniformly over $\ell $ we have,
for $0<\psi <C$ ,%
\begin{equation}
\mathcal{J}_{\ell }(\psi ;4)=O(\ell )\text{, }  \label{PartA}
\end{equation}%
and, for $C<\psi <L\frac{\pi }{2},$%
\begin{equation}
\mathcal{J}_{\ell }(\psi ;4)=\frac{1}{64}\frac{1}{\psi \sin \frac{\psi }{L}}+%
\frac{5}{64}\frac{\cos 4\psi }{\psi \sin \frac{\psi }{L}}-\frac{3}{16}\frac{%
\sin 2\psi }{\psi \sin \frac{\psi }{L}}+O\left(\frac{1}{\psi ^{2}}\frac{1}{\sin
\frac{\psi }{L}}\right)+O\left(\frac{1}{\ell }\frac{1}{\psi \sin \frac{\psi }{L}}\right).  \label{PartB}
\end{equation}
\end{proposition}

The proof of this Proposition is given later in section \ref%
{ProofProposition}; with this result at hand, we can proceed with the proof
of Theorem \ref{maintheorem} as follows.

\subsection{Proof of Theorem \protect\ref{maintheorem}}

\begin{proof}
To establish Theorem \ref{maintheorem}, it is clearly sufficient to show that,
as $\ell \rightarrow \infty ,$%
\begin{equation*}
\corr\left\{ \mathcal{L}_{\ell },\mathcal{M}_{\ell }\right\} =1+O\left(\frac{1}{%
\log \ell }\right)\text{ ,}
\end{equation*}%
and to this end we will prove the equivalent
\begin{equation*}
\cov\left\{ \mathcal{L}_{\ell },\mathcal{M}_{\ell }\right\} =\frac{\log \ell
}{32}+O(1)
\end{equation*}
(cf. \eqref{eq:var nod length log} and \eqref{eq:trispec var log}).
By continuity of the inner product in $L^{2}$ spaces, we need to prove that%
\begin{equation*}
\cov\left\{ \mathcal{L}_{\ell },\mathcal{M}_{\ell }\right\}
=\lim_{\varepsilon \rightarrow 0}\cov\left\{ \mathcal{L}_{\ell ;\varepsilon },%
\mathcal{M}_{\ell }\right\} =\frac{\log \ell }{32}+O(1),
\end{equation*}%
where%
\begin{equation*}
\mathcal{L}_{\ell ;\varepsilon }:=\int_{\mathbb{S}^{2}}\left\Vert \nabla
f_{\ell }(x)\right\Vert \delta _{\varepsilon }(f_{\ell }(x))dx.
\end{equation*}%
Now define the "approximate nodal length"
\begin{equation*}
\Psi _{\varepsilon }(x):=\left\Vert \nabla f_{\ell }(x)\right\Vert \chi
_{\varepsilon }(f_{\ell }(x))\text{ , }\chi _{\varepsilon }(.):=\frac{1}{%
2\varepsilon }\mathbb{I}_{[-\varepsilon ,\varepsilon ]}(.),
\end{equation*}%
where $\mathbb{I}_{A}(.)$ denotes the characteristic function of the set $A$.
The newly defined $\Psi _{\varepsilon }(x)$ is an isotropic random field on $\mathbb{S}^{2}$ admitting the
$L^{2}(\Omega )$ expansion
\begin{equation*}
\Psi _{\varepsilon }(x)=\mathbb{E}\Psi _{\varepsilon }(x)+\sum_{q=2}^{\infty
}\Psi _{\ell ;\varepsilon }(x;q);
\end{equation*}%
moreover, as established in the Appendix, we have the $L^{2}(\Omega )$
convergence
\begin{equation*}
\lim_{\varepsilon \rightarrow 0}\int_{\mathbb{S}^{2}}\Psi _{\varepsilon
}(x)dx=\lim_{\varepsilon \rightarrow 0}\left\{ \left\Vert \nabla f_{\ell
}(x)\right\Vert \chi _{\varepsilon }(f_{\ell }(x))\right\} =\mathcal{L}%
_{\ell }.
\end{equation*}%

Note also that $\Psi _{\varepsilon }(x),H_{4}(f_{\ell }(y))$ are both in $%
L^{2}(\mathbb{S}^{2}\times \Omega )$ and they are isotropic, and thus
\begin{eqnarray*}
\cov\left\{ \mathcal{L}_{\ell ;\varepsilon },\mathcal{M}_{\ell }\right\} &=&-%
\frac{1}{4}\sqrt{\frac{\ell (\ell +1)}{2}}\frac{1}{4!}\cov\left\{ \int_{%
\mathbb{S}^{2}}\Psi _{\varepsilon }(x),\int_{\mathbb{S}^{2}}H_{4}(f_{\ell
}(y))dy\right\} \\
&=&-\frac{1}{4}\sqrt{\frac{\ell (\ell +1)}{2}}\frac{1}{4!}\mathbb{E}\left\{
\int_{\mathbb{S}^{2}}\Psi _{\varepsilon }(x)dx\int_{\mathbb{S}%
^{2}}H_{4}(f_{\ell }(y))dy\right\} \\
&=&-\frac{1}{4}\sqrt{\frac{\ell (\ell +1)}{2}}\frac{1}{4!}\int_{\mathbb{S}%
^{2}}\int_{\mathbb{S}^{2}}\mathbb{E}\left\{ \Psi _{\varepsilon
}(x)H_{4}(f_{\ell }(y))\right\} dxdy \\
&=&-\frac{1}{4}\sqrt{\frac{\ell (\ell +1)}{2}}\frac{1}{4!}\int_{\mathbb{S}%
^{2}}\int_{\mathbb{S}^{2}}\mathbb{E}\left\{ \sum_{q=2}^{\infty }\Psi _{\ell
;\varepsilon }(x;q)H_{4}(f_{\ell }(y))\right\} dxdy \\
&=&\left[ -\frac{1}{4}\sqrt{\frac{\ell (\ell +1)}{2}}\frac{1}{4!}\right]
\times 8\pi ^{2}\int_{0}^{\pi }\mathbb{E}\left\{ \sum_{q=2}^{\infty }\Psi
_{\ell ;\varepsilon }(\overline{x};q)H_{4}(f_{\ell }(y(\theta )))\right\}
\sin \theta d\theta \\
&=&\left[ -\frac{1}{4}\sqrt{\frac{\ell (\ell +1)}{2}}\frac{1}{4!}\right]
\times 8\pi ^{2}\int_{0}^{\pi }\mathbb{E}\left\{ \Psi _{\ell ;\varepsilon }(%
\overline{x};4)H_{4}(f_{\ell }(y(\theta )))\right\} \sin \theta d\theta.
\end{eqnarray*}%
The integrand $\mathbb{E}\left\{ \Psi _{\ell ;\varepsilon }(\overline{x}%
;4)H_{4}(f_{\ell }(y(\theta )))\right\} $ can be computed explicitly and it
is easily seen to be absolutely bounded for fixed $\ell $, uniformly over $%
\varepsilon ,$ see Proposition \ref{Key} below$.$ Hence by the Dominated
Convergence Theorem we may exchange the limit and the integral, and we have that
\begin{eqnarray*}
\cov\left\{ \mathcal{L}_{\ell },\mathcal{M}_{\ell }\right\}
&=&\lim_{\varepsilon \rightarrow 0}\cov\left\{ \mathcal{L}_{\ell
}^{\varepsilon },\mathcal{M}_{\ell }\right\} \\
&=&\left[ -\frac{1}{4}\sqrt{\frac{\ell (\ell +1)}{2}}\frac{1}{4!}\right]
\times 8\pi ^{2}\lim_{\varepsilon \rightarrow 0}\int_{0}^{\pi }\mathbb{E}%
\left\{ \Psi _{\ell ;\varepsilon }(\overline{x};4)H_{4}(f_{\ell }(y(\theta
)))\right\} \sin \theta d\theta \\
&=&\left[ -\frac{1}{4}\sqrt{\frac{\ell (\ell +1)}{2}}\frac{1}{4!}\right]
\times 8\pi ^{2}\int_{0}^{\pi }\lim_{\varepsilon \rightarrow 0}\mathbb{E}%
\left\{ \Psi _{\ell ;\varepsilon }(\overline{x};4)H_{4}(f_{\ell }(y(\theta
)))\right\} \sin \theta d\theta \\
&=&\left[ -\frac{1}{4}\sqrt{\frac{\ell (\ell +1)}{2}}\frac{1}{4!}\right]
\times 8\pi ^{2}\int_{0}^{\pi }\mathbb{E}\left\{ \Psi _{\ell }(\overline{x}%
;4)H_{4}(f_{\ell }(y(\theta )))\right\} \sin \theta d\theta.
\end{eqnarray*}%
We can now rewrite, using (\ref{CrossCorrelation})%
\begin{equation}
\cov\left\{ \mathcal{L}_{\ell },\mathcal{M}_{\ell }\right\} =\int_{0}^{L\pi }%
\mathcal{J}_{\ell }(\psi ;4)\sin \frac{\psi }{L}d\psi.
\label{CovDen}
\end{equation}%
It is now sufficient to notice that
\begin{eqnarray}
\cov\left\{ \mathcal{L}_{\ell },\mathcal{M}_{\ell }\right\}
&=&\int_{0}^{C}\mathcal{J}_{\ell }(\psi ;4)\sin \frac{\psi }{L}d\psi
+2\int_{C}^{L\pi /2}\mathcal{J}_{\ell }(\psi ;4)\sin \frac{\psi }{L}d\psi.
  \label{CovDen2}
\end{eqnarray}
For the first summand in (\ref{CovDen2}) we have easily%
\begin{equation*}
\int_{0}^{C}\mathcal{J}_{\ell }(\psi ;4)\sin \frac{\psi }{L}d\psi \leq \ell
\int_{0}^{C}\sin \frac{\psi }{L}d\psi \leq \frac{\ell }{L}\int_{0}^{C}\psi
d\psi =O(1)\text{ , as }\ell \rightarrow \infty.
\end{equation*}%
For the second sum in (\ref{CovDen2}), using Proposition \ref{Key} and
integrating we obtain%
\begin{equation*}
2\int_{C}^{L\pi /2}\mathcal{J}_{\ell }(\psi ;4)\sin \frac{\psi }{L}d\psi
\end{equation*}%
\begin{eqnarray*}
&=&\frac{1}{L}\int_{C}^{L\pi /2}\frac{1}{\sin ^{2}\frac{\psi }{L}}\left\{
\frac{1}{32}+\frac{5}{32}\cos 4\psi -\frac{3}{8}\sin 2\psi \right\} \sin
\frac{\psi }{L}d\psi \\
&&+O\left(\frac{1}{L}\int_{C}^{L\pi /2}\frac{1}{\psi }\frac{1}{\sin ^{2}\frac{%
\psi }{L}}\sin \frac{\psi }{L}d\psi \right)+O\left(\frac{1}{L}\int_{C}^{L\pi /2}\frac{1%
}{\ell }\frac{1}{\sin ^{2}\frac{\psi }{L}}\sin \frac{\psi }{L}d\psi \right)
\end{eqnarray*}%
\begin{eqnarray*}
&=&\frac{1}{L^{2}}\int_{C}^{L\pi /2}\frac{1}{\psi }\left\{ \frac{1}{32}+%
\frac{5}{32}\cos 4\psi -\frac{3}{8}\sin 2\psi \right\} d\psi \\
&&+O\left(\int_{C}^{L\pi /2}\frac{1}{\psi ^{2}}d\psi \right)+O\left(\frac{1}{L}%
\int_{C}^{L\pi /2}\frac{1}{\ell }\frac{1}{\sin \frac{\psi }{L}}d\psi \right) \\
&=&\frac{\log \ell }{32}+O(1)+O\left(\frac{\log \ell }{\ell }\right),
\end{eqnarray*}%
as claimed.
\end{proof}

\begin{remark}
As mentioned in the Introduction, our main result could equivalently be
stated as

\begin{equation*}
\corr\left\{ \mathcal{L}_{\ell },\proj[\mathcal{L}_{\ell }|4]\right\}
,\corr\left\{ \proj[\mathcal{L}_{\ell }|4],\mathcal{M}_{\ell }\right\}
\rightarrow 1
\end{equation*}%
and thus%
\begin{equation*}
\proj[\widetilde{\mathcal{L}}_{\ell }|4]=-\sqrt{\frac{\ell (\ell +1)}{2}}%
\frac{1}{4\times 24}\int_{\mathbb{S}^{2}}H_{4}(f_{\ell }(x))dx+O(1),
\end{equation*}%
\begin{equation*}
\widetilde{\mathcal{L}}_{\ell }=-\sqrt{\frac{\ell (\ell +1)}{2}}\frac{1}{%
4\times 24}\int_{\mathbb{S}^{2}}H_{4}(f_{\ell }(x))dx+O(1),
\end{equation*}%
both equivalences holding in the $L^{2}(\Omega )$ sense.

\end{remark}

\subsection{Proof of the Central Limit Theorem (Corollary \protect\ref{CLT})}

Recall that $h_{\ell ;4}$ is defined in \eqref{Helle}.
It was shown \cite[Lemma 3.3]{MW2014} that the so-called fourth-order
cumulant of $h_{\ell ;4}$
\begin{equation*}
cum_{4}\left\{ h_{\ell ;4}\right\} :=\mathbb{E}\left[ h_{\ell ;4}^{4}\right]
-3\left\{ \mathbb{E}\left[ h_{\ell ;4}^{4}\right] \right\} ^{2}
\end{equation*}
satisfies $cum_{4}\{h_{\ell ;4}\}\approx \ell ^{-2},$ i.e. the ratio between
the left and right-hand sides is bounded above and below by finite, strictly
positive constants. Taking into account the normalizing factors, it
means that $\widetilde{\mathcal{M}}_{\ell }$ satisfies%
\begin{equation*}
cum_{4}\left\{ \widetilde{\mathcal{M}}_{\ell }\right\} =\mathbb{E}\left[
\widetilde{\mathcal{M}}_{\ell }^{4}\right] -3=O\left(\frac{1}{\log \ell }\right),
\end{equation*}%
where we have exploited the fact that the sequence $\left\{ \widetilde{%
\mathcal{M}}_{\ell }\right\} $ has unit variance. Let us now recall the
so-called \emph{Stein-Malliavin bound} by Nourdin-Peccati, stating that for
a standardized random variable $F$ which belong to the $q-$th order Wiener
chaos $\mathcal{H}_{q}$ we have the bound (see \cite[Theorem 5.2.6]{noupebook})
\begin{equation*}
d_{W}(\widetilde{\mathcal{M}}_{\ell },\mathcal{N}(0,1))\leq
\sqrt{\frac{2q-2}{3\pi q}\left\{ \mathbb{E}F^{4}-3\right\} }.
\end{equation*}

Now the sequence $\left\{ \widetilde{\mathcal{M}}_{\ell }\right\} $ is
indeed standardized and belongs to the Wiener chaos for $q=4,$ so that we
have
\begin{equation*}
d_{W}(\widetilde{\mathcal{M}}_{\ell },\mathcal{N}(0,1))\leq \sqrt{\frac{1}{%
2\pi }\left\{ \mathbb{E}\widetilde{\mathcal{M}}_{\ell }^{4}-3\right\} }=O\left(\frac{1}{\sqrt{\log \ell }}\right).
\end{equation*}%
As a simple application of the triangle inequality for $d_{W}$ (see \cite[Appendix C]{noupebook}), for $\ell \rightarrow \infty $%
\begin{equation*}
d_{W}(\widetilde{\mathcal{L}}_{\ell },\mathcal{N}(0,1))\leq d_{W}(\widetilde{%
\mathcal{M}}_{\ell },\mathcal{N}(0,1))+\sqrt{\mathbb{E}\left[ \widetilde{%
\mathcal{L}}_{\ell }-\widetilde{\mathcal{M}}_{\ell }\right] ^{2}}=O\left(\frac{1}{%
\sqrt{\log \ell }}\right),
\end{equation*}%
and the statement of Corollary \ref{CLT} follows.

\section{\label{ProofProposition} Proof of Proposition \protect\ref{Key}}

\begin{proof}
It is convenient to introduce some further notation,
recalling (\ref{fourthchaos}) and writing%
\begin{equation*}
\Psi _{\ell }(x;4)=A_{\ell }(x)+B_{\ell }(x)+C_{\ell }(x)+D_{\ell
}(x)+E_{\ell }(x)+F_{\ell }(x),
\end{equation*}%
where%
\begin{equation}
\sqrt{\frac{\ell (\ell +1)}{2}}\frac{3}{2}\frac{1}{4!}H_{4}(f_{\ell
}(x))=:A_{\ell }(x),  \label{DefA}
\end{equation}%
\begin{equation}
-\sqrt{\frac{\ell (\ell +1)}{2}}\frac{1}{4}\frac{1}{2!2!}H_{2}(f_{\ell
})H_{2}\left(\frac{\partial _{1;x}f_{\ell }(x)}{\sqrt{\ell (\ell +1)/2}}\right)=:B_{\ell }(x), \label{DefB}
\end{equation}%
\begin{equation}
-\sqrt{\frac{\ell (\ell +1)}{2}}\frac{3}{16}\frac{1}{4!}H_{4}\left(\frac{\partial
_{1;x}f_{\ell }(x)}{\sqrt{\ell (\ell +1)/2}}\right)=:C_{\ell }(x),
\label{DefC}
\end{equation}%
\begin{equation}
+\frac{3}{2}\frac{1}{2!2!}H_{2}\left(\frac{\partial _{1;x}f_{\ell }(x)}{\sqrt{\ell
(\ell +1)/2}}\right)H_{2}\left(\frac{\partial _{2;x}f_{\ell }(x)}{\sqrt{\ell (\ell +1)/2%
}}\right)=:D_{\ell }(x),  \label{DefD}
\end{equation}%
\begin{equation}
-\frac{1}{4}\frac{1}{2!2!}H_{2}(f_{\ell }(x))H_{2}\left(\frac{\partial
_{2;x}f_{\ell }(x)}{\sqrt{\ell (\ell +1)/2}}\right)=:E_{\ell }(x),
\label{DefE}
\end{equation}%
\begin{equation}
-\frac{3}{16}\frac{1}{4!}H_{4}\left(\frac{\partial _{2;x}f_{\ell }(x)}{\sqrt{\ell
(\ell +1)/2}}\right)=:F_{\ell }(x),  \label{DefF}
\end{equation}%
and also
\begin{eqnarray}
\mathcal{M}_{\ell } &:&=-\frac{1}{4}\sqrt{\frac{\ell (\ell +1)}{2}}\frac{1}{%
4!}\int_{\mathbb{S}^{2}}H_{4}(f_{\ell }(x))dx=\int_{\mathbb{S}^{2}}M_{\ell
}(x)dx,  \notag \\
\text{ }M_{\ell }(x) &:&=-\frac{1}{4}\sqrt{\frac{\ell (\ell +1)}{2}}\frac{1}{%
4!}H_{4}(f_{\ell }(x)).  \label{DefM}
\end{eqnarray}%

For the computations to follow, recall that we focus on $\overline{x}=(0,0)$
(the "North Pole") and $y(\theta )=(\theta ,0)$ (the "Greenwich meridian").
By repeated application of the well-known Diagram Formula (see e.g.
\cite[subsection 4.3.1]{MaPeCUP}), we have
\begin{equation*}
\mathbb{E}\left[ H_{2}\left(\frac{\partial _{1;x}f_{\ell }(\overline{x})}{\sqrt{%
\ell (\ell +1)/2}}\right)H_{2}\left(\frac{\partial _{2;x}f_{\ell }(\overline{x})}{\sqrt{%
\ell (\ell +1)/2}}\right)H_{4}(f_{\ell }(y(\theta )))\right]
\end{equation*}%
\begin{equation*}
=4!\frac{4}{\ell ^{2}(\ell +1)^{2}}\left\{ \mathbb{E}\left[ \partial
_{1;x}f_{\ell }(\overline{x})f_{\ell }(y(\theta ))\right] \right\}
^{2}\left\{ \mathbb{E}\left[ \partial _{2;x}f_{\ell }(\overline{x})f_{\ell
}(y(\theta ))\right] \right\} ^{2}=0,
\end{equation*}%
in view of (\ref{cov02}); likewise
\begin{equation*}
\mathbb{E}\left[ H_{2}\left(f_{\ell }(\overline{x})\right)H_{2}\left(\frac{\partial
_{2;x}f_{\ell }(\overline{x})}{\sqrt{\ell (\ell +1)/2}}\right)H_{4}(f_{\ell
}(y(\theta )))\right]
\end{equation*}%
\begin{equation*}
=4!\frac{2}{\ell (\ell +1)}\left\{ \mathbb{E}\left[ f_{\ell }(\overline{x}%
)f_{\ell }(y(\theta ))\right] \right\} ^{2}\left\{ \mathbb{E}\left[ \partial
_{2;x}f_{\ell }(\overline{x})f_{\ell }(y(\theta ))\right] \right\} ^{2}=0%
\text{ ,}
\end{equation*}%
and%
\begin{equation*}
\mathbb{E}\left[ H_{4}\left(\frac{\partial _{2;x}f_{\ell }(x)}{\sqrt{\ell (\ell
+1)/2}}\right)H_{4}(f_{\ell }(y(\theta )))\right]
\end{equation*}%
\begin{equation*}
=4!\frac{4}{\ell ^{2}(\ell +1)^{2}}\left\{ \mathbb{E}\left[ \partial
_{2;x}f_{\ell }(\overline{x})f_{\ell }(y(\theta ))\right] \right\} ^{4}=0.
\end{equation*}%

As a consequence, using definitions (\ref{DefD}), (\ref{DefE}), (\ref{DefF})
and (\ref{DefM}) we have that%
\begin{equation*}
\mathbb{E}\left[ D_{\ell }(\overline{x})M_{\ell }(y(\theta )\right] =\mathbb{%
E}\left[ E_{\ell }(\overline{x})M_{\ell }(y(\theta )\right] =\mathbb{E}\left[
F_{\ell }(\overline{x})M_{\ell }(y(\theta )\right] \equiv 0
\end{equation*}%
for all $\theta \in \lbrack 0,\pi ]$. In the sequel, it is sufficient
to focus on $A_{\ell }(.),B_{\ell }(.)$ and $C_{\ell }(.).$ The proof of \ref%
{PartA} is rather straightforward; indeed, as a simple application of the
Cauchy-Schwartz inequality, we have that%
\begin{equation*}
\mathbb{E}\left[ H_{4}(f_{\ell }(\overline{x}))H_{4}(f_{\ell }(y(\theta )))%
\right] \leq \mathbb{E}\left[ \left\{ H_{4}(f_{\ell }(\overline{x}))\right\}
^{2}\right] =24
\end{equation*}%
\begin{equation*}
\mathbb{E}\left[ H_{2}(f_{\ell }(\overline{x})H_{2}\left(\frac{\partial
_{1;x}f_{\ell }(\overline{x})}{\sqrt{\ell (\ell +1)/2}}\right)H_{4}(f_{\ell
}(y(\theta )))\right]
\end{equation*}%
\begin{equation*}
\leq \sqrt{\mathbb{E}\left[ \left\{ H_{2}(f_{\ell }(\overline{x}))H_{2}\left(\frac{%
\partial _{1;x}f_{\ell }(\overline{x})}{\sqrt{\ell (\ell +1)/2}}\right)\right\}
^{2}\right] \mathbb{E}\left[ \left\{ H_{4}(f_{\ell }(\overline{x}))\right\}
^{2}\right] }=24
\end{equation*}%
and analogously%
\begin{equation*}
\mathbb{E}\left[ H_{4}\left(\frac{\partial _{1;x}f_{\ell }(\overline{x})}{\sqrt{%
\ell (\ell +1)/2}}\right)H_{4}(f_{\ell }(y(\theta )))\right]
\end{equation*}%
\begin{equation*}
\leq \sqrt{\mathbb{E}\left[ \left\{ H_{4}\left(\frac{\partial _{1;x}f_{\ell }(%
\overline{x})}{\sqrt{\ell (\ell +1)/2}}\right)\right\} ^{2}\right] \mathbb{E}\left[
\left\{ H_{4}(f_{\ell }(\overline{x}))\right\} ^{2}\right] }=24.
\end{equation*}%
It then follows that
\begin{equation*}
\left\vert \mathcal{J}_{\ell }(\psi ;4)\right\vert =8\pi ^{2}\left\vert
\mathbb{E}\left[ \left\{ A_{\ell }(\overline{x})+B(\overline{x})+C_{\ell }(%
\overline{x})\right\} M_{\ell }(y(\frac{\psi }{L}))\right] \right\vert
\end{equation*}%
\begin{equation*}
\leq \frac{8\pi ^{2}}{L}\left\{ \left\vert \mathbb{E}\left[ A_{\ell }(%
\overline{x})M_{\ell }(y(\frac{\psi }{L}))\right] \right\vert +\left\vert
\mathbb{E}\left[ B(\overline{x})M_{\ell }(y(\frac{\psi }{L}))\right]
\right\vert +\left\vert \mathbb{E}\left[ C_{\ell }(\overline{x})M_{\ell }(y(%
\frac{\psi }{L}))\right] \right\vert \right\}
\end{equation*}%
\begin{equation*}
\leq 24\frac{\ell (\ell +1)}{2L}8\pi ^{2}\left\{ \frac{3}{2}\frac{1}{4!}%
\frac{1}{4}\frac{1}{4!}+\frac{1}{4}\frac{1}{2!2!}\frac{1}{4}\frac{1}{4!}+%
\frac{3}{16}\frac{1}{4!}\frac{1}{4}\frac{1}{4!}\right\} =O(\ell ),
\end{equation*}%
as claimed.

We now turn to proving \ref{PartB}. Using \eqref{cov00}, \eqref{cov01} and the Diagram Formula we can write explicitly
\begin{eqnarray}
\mathbb{E}\left[ A_{\ell }(\overline{x})M_{\ell }(y(\theta ))\right] &=&-%
\frac{\ell (\ell +1)}{2}\frac{3}{2}\frac{1}{4!}\frac{1}{4}\frac{1}{4!}\times
4!P_{\ell }^{4}(\cos \theta )  \notag \\
&=&-\frac{\ell (\ell +1)}{2}\frac{1}{64}P_{\ell }^{4}(\cos \theta )\text{ ,}
\label{Aemme}
\end{eqnarray}
\begin{eqnarray}
\mathbb{E}\left[ B_{\ell }(\overline{x})M_{\ell }(y(\theta ))\right] &=&%
\frac{\ell (\ell +1)}{2}\frac{1}{4}\frac{1}{2!2!}\frac{1}{4}\frac{1}{4!}%
\times 4!\frac{2}{\ell (\ell +1)}P_{\ell }^{2}(\cos \theta )\left\{ P_{\ell
}^{\prime }(\cos \theta )\right\} ^{2}  \notag \\
&=&\frac{\ell (\ell +1)}{2}\frac{1}{64}\frac{2}{\ell (\ell +1)}P_{\ell
}^{2}(\cos \theta )\left\{ P_{\ell }^{\prime }(\cos \theta )\right\} ^{2},  \label{Bemme}
\end{eqnarray}
\begin{eqnarray}
\mathbb{E}\left[ C_{\ell }(\overline{x})M_{\ell }(y(\theta ))\right] &=&%
\frac{\ell (\ell +1)}{2}\frac{3}{16}\frac{1}{4!}\frac{1}{4}\frac{1}{4!}%
\times \frac{4}{\ell ^{2}(\ell +1)^{2}}4!\left\{ P_{\ell }^{\prime }(\cos
\theta )\sin \theta \right\} ^{4}  \notag \\
&=&\frac{\ell (\ell +1)}{2}\frac{3}{16}\frac{1}{4!}\frac{1}{\ell ^{2}(\ell
+1)^{2}}\left\{ P_{\ell }^{\prime }(\cos \theta )\sin \theta \right\} ^{4}.  \label{Cemme}
\end{eqnarray}
Now recall that (see section \ref{Appendix})%
\begin{equation*}
P_{\ell }(\cos \frac{\psi }{L})=\sqrt{\frac{2}{\pi \ell \sin \frac{\psi }{L}}%
}\left(\sin \left(\psi +\frac{\pi }{4}\right)+O\left(\frac{1}{\psi }\right)\right)
\end{equation*}%
\begin{equation*}
P_{\ell }^{\prime }\left(\cos \frac{\psi }{L}\right)=\sqrt{\frac{2}{\pi \ell \sin ^{3}%
\frac{\psi }{L}}}\left(\ell \sin \left(\psi -\frac{\pi }{4}\right)+O(1)\right).
\end{equation*}%
Let us also mention the following standard trigonometric identities that are
used repeatedly in our arguments:%
\begin{equation*}
\left\{ \sin (\psi -\frac{\pi }{4})\right\} ^{4}=\left\{ \frac{\sqrt{2}}{2}%
\sin \psi -\frac{\sqrt{2}}{2}\cos \psi \right\} ^{4}=\frac{3}{8}-\frac{1}{8}%
\cos 4\psi -\frac{1}{2}\sin 2\psi,
\end{equation*}%
\begin{equation*}
\left\{ \sin (\psi +\frac{\pi }{4})\right\} ^{4}=\frac{3}{8}-\frac{1}{8}\cos
4\psi +\frac{1}{2}\sin 2\psi,
\end{equation*}%
and%
\begin{equation*}
(1+\sin 2\psi )(1-\sin 2\psi )=\frac{1+\cos 4\psi }{2}\text{ .}
\end{equation*}%

Substituting the latter expressions into (\ref{Aemme}), we obtain that%
\begin{equation*}
8\pi ^{2}\mathbb{E}\left[ A_{\ell }(\overline{x})M_{\ell }(y(\theta ))\right]
=-\frac{\ell (\ell +1)}{2}\frac{1}{8}\pi ^{2}P_{\ell }^{4}(\cos \theta )
\end{equation*}%
\begin{eqnarray*}
&=&-\frac{\ell (\ell +1)}{2}\frac{1}{8}\pi ^{2}\left[ \sqrt{\frac{2}{\pi
\ell \sin \frac{\psi }{L}}}(\sin (\psi +\frac{\pi }{4})+O\left(\frac{1}{\psi }\right))%
\right] ^{4} \\
&=&-\frac{\ell (\ell +1)}{2}\frac{1}{8}\pi ^{2}\frac{2^{2}}{\pi ^{2}\ell
^{2}\sin ^{2}\frac{\psi }{L}}\left\{ \sin (\psi +\frac{\pi }{4})\right\}
^{4}+O\left(\frac{1}{\psi }\frac{1}{\sin ^{2}\frac{\psi }{L}}\right) \\
&=&-\frac{\ell (\ell +1)}{2}\frac{1}{2\ell ^{2}\sin ^{2}\frac{\psi }{L}}%
\left\{ \frac{3}{8}-\frac{1}{8}\cos 4\psi +\frac{1}{2}\sin 2\psi \right\} +O\left(
\frac{1}{\psi }\frac{1}{\sin ^{2}\frac{\psi }{L}}\right) \\
&=&-\frac{1}{4\sin ^{2}\frac{\psi }{L}}\left\{ \frac{3}{8}-\frac{1}{8}\cos
4\psi +\frac{1}{2}\sin 2\psi \right\} +O\left(\frac{1}{\psi }\frac{1}{\sin ^{2}%
\frac{\psi }{L}}\right)+O\left(\frac{1}{\ell }\frac{1}{\sin ^{2}\frac{\psi }{L}}\right).
\end{eqnarray*}

Likewise for (\ref{Bemme})%
\begin{equation*}
8\pi ^{2}\mathbb{E}\left[ B_{\ell }(\overline{x})M_{\ell }(y(\theta ))\right]
\end{equation*}%
\begin{equation*}
=\frac{\ell (\ell +1)}{2}\frac{1}{8}\pi ^{2}\frac{2}{\ell (\ell +1)}P_{\ell
}^{2}(\cos \theta )\left\{ P_{\ell }^{\prime }(\cos \theta )\sin \theta
\right\} ^{2}
\end{equation*}%
\begin{equation*}
=\frac{1}{8}\pi ^{2}\left[ \sqrt{\frac{2}{\pi \ell \sin \frac{\psi }{L}}}%
\left(\sin (\psi +\frac{\pi }{4})+O\left(\frac{1}{\psi }\right)\right)\right] ^{2}\left[ \sqrt{%
\frac{2}{\pi \ell \sin ^{3}\frac{\psi }{L}}}\left(\ell \sin \left(\psi -\frac{\pi }{4}%
\right)+O(1)\right)\sin \frac{\psi }{L}\right] ^{2}
\end{equation*}%
\begin{equation*}
=\frac{1}{8}\pi ^{2}\frac{2}{\pi \ell \sin \frac{\psi }{L}}\sin ^{2}(\psi +%
\frac{\pi }{4})\frac{2}{\pi \ell \sin \frac{\psi }{L}}\ell ^{2}\sin
^{2}\left(\psi -\frac{\pi }{4}\right)+O\left(\frac{1}{\ell \sin ^{2}\frac{\psi }{L}}\right)
\end{equation*}%
\begin{equation*}
=\frac{1}{2}\frac{1}{\sin \frac{\psi }{L}}\sin ^{2}(\psi +\frac{\pi }{4})%
\frac{1}{\sin \frac{\psi }{L}}\sin ^{2}\left(\psi -\frac{\pi }{4}\right)+O\left(\frac{1}{%
\ell \sin ^{2}\frac{\psi }{L}}\right)
\end{equation*}%
\begin{equation*}
=\frac{1}{2}\frac{1}{\sin \frac{\psi }{L}}\left\{ \frac{\sqrt{2}}{2}\sin
\psi -\frac{\sqrt{2}}{2}\cos \psi \right\} ^{2}\frac{1}{\sin \frac{\psi }{L}}%
\left\{ \frac{\sqrt{2}}{2}\sin \psi -\frac{\sqrt{2}}{2}\cos \psi \right\}
^{2}+O\left(\frac{1}{\ell \sin ^{2}\frac{\psi }{L}}\right)
\end{equation*}%
\begin{equation*}
=\frac{1}{2}\frac{1}{\sin \frac{\psi }{L}}\frac{1}{4}\left\{ \sin \psi +\cos
\psi \right\} ^{2}\frac{1}{\sin \frac{\psi }{L}}\left\{ \sin \psi -\cos \psi
\right\} ^{2}+O\left(\frac{1}{\ell \sin ^{2}\frac{\psi }{L}}\right)
\end{equation*}%
\begin{equation*}
=\frac{1}{2}\frac{1}{\sin ^{2}\frac{\psi }{L}}\frac{1}{4}\left\{ \sin
^{2}\psi -\cos ^{2}\psi \right\} ^{2}+O\left(\frac{1}{\ell \sin ^{2}\frac{\psi }{L%
}}\right)=\frac{1}{2}\frac{1}{\sin ^{2}\frac{\psi }{L}}\frac{1}{4}\left\{ 1-2\cos
^{2}\psi \right\} ^{2}+O\left(\frac{1}{\ell \sin ^{2}\frac{\psi }{L}}\right)
\end{equation*}%
\begin{equation*}
=\frac{1}{2}\frac{1}{\sin ^{2}\frac{\psi }{L}}\frac{1+\cos 4\psi }{8}+O\left(
\frac{1}{\ell \sin ^{2}\frac{\psi }{L}}\right).
\end{equation*}%
Finally, for (\ref{Cemme})%
\begin{equation*}
8\pi ^{2}\mathbb{E}\left[ C_{\ell }(\overline{x})M_{\ell }(y(\theta ))\right]
\end{equation*}%
\begin{equation*}
=\frac{\ell (\ell +1)}{2}\frac{3}{2}\pi ^{2}\frac{1}{4!}\frac{1}{\ell
^{2}(\ell +1)^{2}}\left\{ P_{\ell }^{\prime }(\cos \theta )\sin \theta
\right\} ^{4}
\end{equation*}%
\begin{equation*}
=\frac{3}{4}\frac{1}{4!}\pi ^{2}\frac{1}{\ell (\ell +1)}\left\{ \sqrt{\frac{2%
}{\pi \ell \sin ^{3}\frac{\psi }{L}}}\left(\ell \sin \left(\psi -\frac{\pi }{4}%
\right)+O(1)\right)\sin \frac{\psi }{L}\right\} ^{4}
\end{equation*}%
\begin{equation*}
=\frac{3}{4}\frac{1}{4!}\pi ^{2}\frac{1}{\ell (\ell +1)}\frac{2^{2}}{\pi
^{2}\ell ^{2}\sin ^{2}\frac{\psi }{L}}\ell ^{4}\sin ^{4}\left(\psi -\frac{\pi }{4}%
\right)+O\left(\frac{1}{\ell \sin ^{2}\frac{\psi }{L}}\right)
\end{equation*}%
\begin{eqnarray*}
&=&\frac{1}{8}\frac{1}{\sin ^{2}\frac{\psi }{L}}\left[ \frac{3}{8}-\frac{1}{8%
}\cos 4\psi -\frac{1}{2}\sin 2\psi \right] +O\left(\frac{1}{\ell \sin ^{2}\frac{%
\psi }{L}}\right) \\
&=&\frac{1}{8}\frac{1}{\sin ^{2}\frac{\psi }{L}}\left[ \frac{3}{8}-\frac{1}{8%
}\cos 4\psi -\frac{1}{2}\sin 2\psi \right] +O\left(\frac{1}{\ell \sin ^{2}\frac{%
\psi }{L}}\right)\text{ .}
\end{eqnarray*}%
Thus, summing (\ref{Aemme}), (\ref{Bemme}) and (\ref{Cemme}) we obtain
\begin{eqnarray*}
\mathcal{J}_{\ell }(\psi ;4) &=&-\frac{1}{4L\sin ^{2}\frac{\psi }{L}}\left\{
\frac{3}{8}-\frac{1}{8}\cos 4\psi +\frac{1}{2}\sin 2\psi \right\} +O\left(\frac{1%
}{\psi }\frac{1}{\sin ^{2}\frac{\psi }{L}}\right)+O\left(\frac{1}{\ell }\frac{1}{\sin
^{2}\frac{\psi }{L}}\right) \\
&&+\frac{1}{2}\frac{1}{L\sin ^{2}\frac{\psi }{L}}\frac{1+\cos 4\psi }{8}+O\left(
\frac{1}{\ell \sin ^{2}\frac{\psi }{L}}\right) \\
&&+\frac{1}{8}\frac{1}{L\sin ^{2}\frac{\psi }{L}}\left[ \frac{3}{8}-\frac{1}{%
8}\cos 4\psi -\frac{1}{2}\sin 2\psi \right] +O\left(\frac{1}{\ell \sin ^{2}\frac{%
\psi }{L}}\right),
\end{eqnarray*}%
\begin{equation*}
=\frac{1}{64}\frac{1}{L\sin ^{2}\frac{\psi }{L}}+\frac{5}{64}\frac{\cos
4\psi }{L\sin ^{2}\frac{\psi }{L}}-\frac{3}{16}\frac{\sin 2\psi }{L\sin ^{2}%
\frac{\psi }{L}}+O\left(\frac{1}{\psi }\frac{1}{L\sin ^{2}\frac{\psi }{L}}\right)+O\left(\frac{1}{\ell }\frac{1}{L\sin ^{2}\frac{\psi }{L}}\right),
\end{equation*}%
\begin{equation}
=\frac{1}{64}\frac{1}{\psi \sin \frac{\psi }{L}}+\frac{5}{64}\frac{\cos
4\psi }{\psi \sin \frac{\psi }{L}}-\frac{3}{16}\frac{\sin 2\psi }{\psi \sin
\frac{\psi }{L}}+O\left(\frac{1}{\psi ^{2}}\frac{1}{\sin \frac{\psi }{L}}\right)+O\left(\frac{1}{\ell }\frac{1}{\psi \sin \frac{\psi }{L}}\right), \label{KeyProp}
\end{equation}%
as claimed.
\end{proof}

\begin{remark}
The variance of the spherical nodal length is written \cite[Proposition 2.7]{Wig} as
\begin{equation*}
\var\left\{ \mathcal{L}_{\ell }\right\} =\int_{0}^{L\pi }4\pi ^{2}\frac{\ell
(\ell +1)}{L}\left\{ \mathcal{K}_{\ell }(\psi )-\frac{1}{4}\right\} \sin (%
\frac{\psi }{L})d\psi.
\end{equation*}
Here $\mathcal{K}_{\ell }(\cdot )$ represents the
two-point correlation function of the nodal length, defined as
\begin{equation*}
\mathcal{K}_{\ell }(\psi )=\frac{1}{2\pi \sqrt{1-P_{\ell }^{2}(\cos \frac{%
\psi }{L})}}\mathbb{E}\left[ \left. \left\Vert \nabla f_{\ell
}(N)\right\Vert  \cdot \left\Vert \nabla f_{\ell }(\frac{\psi }{L})\right\Vert
\right\vert f_{\ell }(N)=f_{\ell }\left(\frac{\psi }{L}\right)=0\right].
\end{equation*}%
It was shown ~\cite{Wig} that one has%
\begin{equation*}
\mathcal{K}_{\ell }(\psi )-\frac{1}{4}=\frac{1}{2}\frac{\sin 2\psi }{\pi
\ell \sin \frac{\psi }{L}}+\frac{9}{32}\frac{\cos 2\psi }{\pi \ell \psi \sin
\frac{\psi }{L}}+\frac{1}{256}\frac{1}{\pi ^{2}\ell \psi \sin \frac{\psi }{L}}
\end{equation*}%
\begin{equation*}
+\frac{\frac{27}{64}\sin 2\psi -\frac{75}{256}\cos 4\psi }{\pi ^{2}\ell \psi
\sin \frac{\psi }{L}}+O\left(\frac{1}{\psi ^{3}}+\frac{1}{\ell \psi }\right).
\end{equation*}%
To compare this result with those in the present paper, let us note that
\begin{equation*}
4\pi ^{2}\frac{\ell (\ell +1)}{L}\left\{ \mathcal{K}_{\ell }(\psi )-\frac{1}{%
4}\right\} =\frac{1}{64}\cdot \frac{1}{\psi \sin \frac{\psi }{L}}+\text{%
oscillatory or lower order terms ,}
\end{equation*}%
in perfect analogy with the two-point
cross-correlation function \eqref{KeyProp}, used to compute the covariance
\begin{equation*}
\cov\left\{ \mathcal{L}_{\ell },\mathcal{M}_{\ell }\right\} =\int_{0}^{L\pi }%
\mathcal{J}_{\ell }(\psi ;4)\sin \left(\frac{\psi }{L}\right)d\psi,
\end{equation*}
satisfying
\begin{equation*}
\mathcal{J}_{\ell }(\psi ;4)=\frac{1}{64}\frac{1}{\psi \sin \frac{\psi }{L}}+%
\text{oscillatory or lower order terms.}
\end{equation*}
\end{remark}

\appendix

\label{Appendix}

\section{Some Background Material}

For completeness, in this Appendix we record some basic facts about
covariances of random spherical harmonics and their derivatives; all the
expressions to follow are rather standard and have been repeatedly exploited
in the literature.
Let us first recall that for arbitrary coordinates $x=(\theta _{x},\varphi
_{x}),$ $y=(\theta _{y},\varphi _{y})$ we have%
\begin{equation*}
\left\langle x,y\right\rangle =\cos \theta _{x}\cos \theta _{y}+\sin \theta
_{x}\sin \theta _{y}\cos (\varphi _{x}-\varphi _{y})
\end{equation*}%
It is then elementary to show that%
\begin{equation*}
\mathbb{E}\left[ f_{\ell }(x)\partial _{1;y}f_{\ell }(y)\right] =P_{\ell
}^{\prime }(\left\langle x,y\right\rangle )\left\{ -\cos \theta _{x}\sin
\theta _{y}+\sin \theta _{x}\cos \theta _{y}\cos (\varphi _{x}-\varphi
_{y})\right\} \text{ ,}
\end{equation*}%
\begin{equation*}
\mathbb{E}\left[ f_{\ell }(x)\partial _{2;y}f_{\ell }(y)\right] =P_{\ell
}^{\prime }(\left\langle x,y\right\rangle )\sin \theta _{x}\sin (\varphi
_{x}-\varphi _{y})\text{ ,}
\end{equation*}%
\begin{equation*}
\mathbb{E}\left[ \partial _{1;x}f_{\ell }(x)\partial _{1;y}f_{\ell }(y)%
\right]
\end{equation*}%
\begin{eqnarray*}
&=&P_{\ell }^{\prime \prime }(\left\langle x,y\right\rangle )\left\{ -\cos
\theta _{x}\sin \theta _{y}+\sin \theta _{x}\cos \theta _{y}\cos (\varphi
_{x}-\varphi _{y})\right\} \left\{ -\sin \theta _{x}\cos \theta _{y}+\cos
\theta _{x}\sin \theta _{y}\cos (\varphi _{x}-\varphi _{y})\right\} \\
&&+P_{\ell }^{\prime }(\left\langle x,y\right\rangle )\left\{ \sin \theta
_{x}\sin \theta _{y}+\cos \theta _{x}\cos \theta _{y}\cos (\varphi
_{x}-\varphi _{y})\right\},
\end{eqnarray*}%
\begin{eqnarray*}
\mathbb{E}\left[ \partial _{1;x}f_{\ell }(x)\partial _{2;y}f_{\ell }(y)%
\right] &=&-P_{\ell }^{\prime \prime }(\left\langle x,y\right\rangle
)\left\{ \sin \theta _{x}\cos \theta _{y}+\cos \theta _{x}\sin \theta
_{y}\cos (\varphi _{x}-\varphi _{y})\right\} \sin \theta _{x}\sin (\varphi
_{x}-\varphi _{y}) \\
&&-P_{\ell }^{\prime }(\left\langle x,y\right\rangle )\cos \theta _{x}\sin
(\varphi _{x}-\varphi _{y})\text{ ,}
\end{eqnarray*}%
\begin{equation*}
\mathbb{E}\left[ \partial _{2;x}f_{\ell }(x)\partial _{2;y}f_{\ell }(y)%
\right] =-P_{\ell }^{\prime \prime }(\left\langle x,y\right\rangle )\sin
\theta _{x}\sin \theta _{y}\sin ^{2}(\varphi _{x}-\varphi _{y})+P_{\ell
}^{\prime }(\left\langle x,y\right\rangle )\cos (\varphi _{x}-\varphi _{y})~.
\end{equation*}

In particular, the result we exploited several times in this paper are
obtained setting $x=(0,0),$ $y=(\theta ,0)$:
\begin{equation}
\mathbb{E}\left[ f_{\ell }(x)f_{\ell }(y)\right] =P_{\ell }(\cos \theta ), \label{cov00}
\end{equation}%
\begin{equation}
\mathbb{E}\left[ f_{\ell }(x)\partial _{1;y}f_{\ell }(y)\right] =-P_{\ell
}^{\prime }(\cos \theta )\sin \theta,  \label{cov01}
\end{equation}%
\begin{equation}
\mathbb{E}\left[ f_{\ell }(x)\partial _{2;y}f_{\ell }(y)\right] =\mathbb{E}%
\left[ f_{\ell }(y)\partial _{2;y}f_{\ell }(x)\right] =0,
\label{cov02}
\end{equation}%
\begin{equation}
\mathbb{E}\left[ \partial _{1;x}f_{\ell }(x)\partial _{1;y}f_{\ell }(y)%
\right] =P_{\ell }^{\prime }(\cos \theta )\cos \theta -P_{\ell }^{\prime
\prime }(\cos \theta )\sin ^{2}\theta \text{ ,}  \label{cov11}
\end{equation}%
\begin{equation}
\mathbb{E}\left[ \partial _{1;x}f_{\ell }(x)\partial _{2;y}f_{\ell }(y)%
\right] =\mathbb{E}\left[ \partial _{1;x}f_{\ell }(y)\partial _{2;y}f_{\ell
}(x)\right] =0\text{ ,}  \label{cov12}
\end{equation}%
\begin{equation}
\mathbb{E}\left[ \partial _{2;x}f_{\ell }(x)\partial _{2;y}f_{\ell }(y)%
\right] =P_{\ell }^{\prime }(\cos \theta ).  \label{cov22}
\end{equation}%
On the other hand, the following very useful expansions are proved \cite[LemmaB.3]{Wig},
and hold uniformly for $C<\psi <L\frac{\pi }{2}$ (recall
that $L:=\ell +\frac{1}{2}):$%
\begin{equation*}
P_{\ell }(\cos \frac{\psi }{L})=\sqrt{\frac{2}{\pi \ell \sin \frac{\psi }{L}}%
}\left(\sin \left(\psi +\frac{\pi }{4}\right)+O\left(\frac{1}{\psi }\right)\right),
\end{equation*}%
\begin{equation*}
P_{\ell }^{\prime }(\cos \frac{\psi }{L})=\sqrt{\frac{2}{\pi \ell \sin ^{3}%
\frac{\psi }{L}}}(\ell \sin \left(\psi -\frac{\pi }{4}\right)+O(1)),
\end{equation*}%
\begin{equation*}
P_{\ell }^{\prime \prime }(\cos \frac{\psi }{L})=-\frac{\ell ^{2}}{\sin ^{2}%
\frac{\psi }{L}}P_{\ell }(\cos \frac{\psi }{L})+\frac{2}{\sin ^{2}\frac{\psi
}{L}}P_{\ell }^{\prime }\left(\cos \frac{\psi }{L}\right)+O\left(\frac{\ell ^{3}}{\psi ^{5/2}}\right).
\end{equation*}

\section{\emph{The }$L^{2}$\emph{\ approximation}}

We know that the nodal length is defined almost-surely by
\begin{equation*}
\lim_{\varepsilon \rightarrow 0}\int_{\mathbb{S}^{2}}\chi _{\varepsilon
}(f_{\ell }(x))\left\Vert \nabla f_{\ell }(x)\right\Vert dx;
\end{equation*}%
the almost-sure convergence follows from the standard argument (\cite[Lemma 3.1]{RudWig}), as
$\chi _{\varepsilon }(.)=\frac{1}{2\varepsilon }\mathbb{I}_{[-\varepsilon ,\varepsilon ]}(.)$ is integrable
and $f_{\ell }$ is smooth we have, using the co-area formula \cite[p.169]{adlertaylor}
\begin{equation*}
\int_{\mathbb{S}^{2}}\chi _{\varepsilon }(f_{\ell }(x))\left\Vert \nabla
f_{\ell }(x)\right\Vert dx=\int_{\mathbb{R}}\left\{ \int_{f_{\ell
}^{-1}(s)}\chi _{\varepsilon }(f_{\ell }(x))\right\} ds.
\end{equation*}
Since
\begin{equation*}
\chi _{\varepsilon }(f_{\ell }(x))=\left\{
\begin{array}{c}
0\text{ for }x:f_{\ell }(x)>\varepsilon \\
\frac{1}{2\varepsilon }\text{ for }x:f_{\ell }(x)\leq \varepsilon%
\end{array}%
\right.
\end{equation*}%
we obtain%
\begin{equation*}
\int_{\mathbb{R}}\left\{ \int_{f_{\ell }^{-1}(s)}\chi _{\varepsilon
}(f_{\ell }(x))\right\} ds=\frac{1}{2\varepsilon }\int_{-\varepsilon
}^{\varepsilon }\vol\left[ f_{\ell }^{-1}(s)\right] ds\rightarrow \vol\left[
f_{\ell }^{-1}(0)\right] \text{ , as }\varepsilon \rightarrow 0\text{ ,}
\end{equation*}%
since the function $s\rightarrow \vol\left[ f_{\ell }^{-1}(s)\right] $ is
continuous for regular (Morse) functions. We now want to show that the
convergence occurs also in the $L^{2}$ sense; as convergence holds
almost surely, it is sufficient to show that
\begin{equation*}
\lim_{\varepsilon \rightarrow 0}\mathbb{E}\left[ \mathcal{L}_{\ell
;\varepsilon }^{2}\right] =\mathbb{E}\left[ \mathcal{L}_{\ell }^{2}\right].
\end{equation*}%
Note that%
\begin{eqnarray*}
\mathbb{E}\left[ \mathcal{L}_{\ell ;\varepsilon }^{2}\right] &=&\mathbb{E}%
\left[ \left\{ \int_{\mathbb{S}^{2}}\left\{ \chi _{\varepsilon }(f_{\ell
}(x))\left\Vert \nabla f_{\ell }(x)\right\Vert \right\} dx\right\} ^{2}%
\right] \\
&=&\mathbb{E}\left[ \left\{ \int_{\mathbb{R}}\int_{f_{\ell }(x)=u}\chi
_{\varepsilon }(f_{\ell }(x))dxdu\right\} ^{2}\right] \\
&=&\mathbb{E}\left[ \left\{ \int_{\mathbb{R}}\mathcal{L}_{\ell }(u)\chi
_{\varepsilon }(u)du\right\} ^{2}\right].
\end{eqnarray*}

It is easy to see that the application $u\rightarrow \mathbb{E}\left[
\left\{ \mathcal{L}_{\ell }(u)\right\} ^{2}\right] $ is continuous, where%
\begin{eqnarray*}
\mathbb{E}\left[ \mathcal{L}_{\ell }^{2}(u)\right] &=&\int_{\mathbb{S}%
^{2}\times \mathbb{S}^{2}}\mathbb{E}\left[ \left. \left\Vert \nabla f_{\ell
}(x_{1})\right\Vert \left\Vert \nabla f_{\ell }(x_{2})\right\Vert
\right\vert f_{\ell }(x_{1})=u,f_{\ell }(x_{2})=u\right] \phi _{f_{\ell
}(x_{1}),f_{\ell }(x_{2})}(u,u)dx_{1}dx_{2} \\
&=&8\pi ^{2}\int_{0}^{\pi }\text{ }\mathbb{E}\left[ \left. \left\Vert \nabla
f_{\ell }(N)\right\Vert \left\Vert \nabla f_{\ell }(y(\theta ))\right\Vert
\right\vert f_{\ell }(N)=u,f_{\ell }(y(\theta ))=u\right] \phi _{f_{\ell
}(N),f_{\ell }(y(\theta ))}(u,u)\sin \theta d\theta.
\end{eqnarray*}%
To check continuity, it is enough to show that the Dominated Convergence Theorem holds; we first note that%
\begin{eqnarray*}
\phi _{f_{\ell }(N),f_{\ell }(y(\theta ))}(u,u)\sin \theta &\leq &\phi
_{f_{\ell }(N),f_{\ell }(y(\theta ))}(0,0)\sin \theta \\
&=&\frac{1}{2\pi \sqrt{1-P_{\ell }^{2}(\cos \theta )}}\sin \theta =O(1)\text{
,}
\end{eqnarray*}%
uniformly over $\theta $. On the other hand, to evaluate%
\begin{equation*}
\mathbb{E}\left[ \left. \left\Vert \nabla f_{\ell }(x_{1})\right\Vert
\left\Vert \nabla f_{\ell }(x_{2})\right\Vert \right\vert f_{\ell
}(N)=u,f_{\ell }(y(\theta ))=u\right]
\end{equation*}%
we can use Cauchy-Schwartz inequality, and bound
\begin{equation*}
\mathbb{E}\left[ \left. w_{i}^{2}\right\vert f_{\ell }(N)=u,f_{\ell
}(y(\theta ))=u\right] =\var\left[ \left. w_{i}\right\vert f_{\ell
}(N)=u,f_{\ell }(y(\theta ))=u\right] +\left\{ \mathbb{E}\left[ \left.
w_{i}\right\vert f_{\ell }(N)=u,f_{\ell }(y(\theta ))=u\right] \right\} ^{2},
\end{equation*}%
for $i=1,2,3,4,$ where%
\begin{equation*}
\left(
\begin{array}{c}
w_{1} \\
w_{2} \\
w_{3} \\
w_{4}%
\end{array}%
\right) :=\left(
\begin{array}{c}
\nabla f_{\ell }(x_{1}) \\
\nabla f_{\ell }(x_{2})%
\end{array}%
\right).
\end{equation*}%
Note first that, by standard properties of Gaussian conditional distributions%
\begin{equation*}
\var\left[ \left. w_{i}\right\vert f_{\ell }(N)=u,f_{\ell }(y(\theta ))=u%
\right] =\var\left[ \left. w_{i}\right\vert f_{\ell }(N)=0,f_{\ell }(y(\theta
))=0\right],
\end{equation*}%
and the quantities on the right-hand sides have been shown to be uniformly
bounded over $\theta $ in \cite{Wig}. On the other hand, a direct computation
along the same lines as in \cite[Appendix A]{Wig} shows that%
\begin{equation*}
\mathbb{E}\left[ \left.
\begin{array}{c}
w_{1} \\
w_{2} \\
w_{3} \\
w_{4}%
\end{array}%
\right\vert f_{\ell }(N)=u,f_{\ell }(y(\theta ))=u\right] =B_{\ell
}^{T}(\theta )A_{\ell }^{-1}(\theta )\left(
\begin{array}{c}
u \\
u%
\end{array}%
\right) \text{ ,}
\end{equation*}%
where%
\begin{eqnarray*}
B_{\ell }^{T}(\theta ) &=&\left(
\begin{array}{cc}
-P_{\ell }^{\prime }(\cos \theta )\sin \theta & 0 \\
0 & 0 \\
0 & P_{\ell }^{\prime }(\cos \theta )\sin \theta \\
0 & 0%
\end{array}%
\right) \text{ , } \\
A_{\ell }^{-1}(\theta ) &=&\frac{1}{1-P_{\ell }^{2}(\cos \theta )}\left(
\begin{array}{cc}
1 & -P_{\ell }(\cos \theta ) \\
-P_{\ell }(\cos \theta ) & 1%
\end{array}%
\right),
\end{eqnarray*}%
so that the conditional expected value can be written as%
\begin{equation*}
\frac{1}{1-P_{\ell }^{2}(\cos \theta )}\left(
\begin{array}{cc}
-P_{\ell }^{\prime }(\cos \theta )\sin \theta & P_{\ell }^{\prime }(\cos
\theta )P_{\ell }(\cos \theta )\sin \theta \\
0 & 0 \\
-P_{\ell }^{\prime }(\cos \theta )P_{\ell }(\cos \theta )\sin \theta &
P_{\ell }^{\prime }(\cos \theta )\sin \theta \\
0 & 0%
\end{array}%
\right) \left(
\begin{array}{c}
u \\
u%
\end{array}%
\right)
\end{equation*}%
\begin{equation*}
=\frac{1}{1-P_{\ell }^{2}(\cos \theta )}\left(
\begin{array}{c}
uP_{\ell }^{\prime }(\cos \theta )\sin \theta (P_{\ell }(\cos \theta )-1) \\
0 \\
uP_{\ell }^{\prime }(\cos \theta )\sin \theta (1-P_{\ell }(\cos \theta )) \\
0%
\end{array}%
\right)
\end{equation*}%
\begin{equation*}
=\frac{1}{1+P_{\ell }(\cos \theta )}\left(
\begin{array}{c}
-uP_{\ell }^{\prime }(\cos \theta )\sin \theta \\
0 \\
uP_{\ell }^{\prime }(\cos \theta )\sin \theta \\
0%
\end{array}%
\right) \text{ .}
\end{equation*}%

This vector function is immediately seen to be uniformly bounded over $%
\theta ,$ whence the Dominated Convergence Theorem holds. Hence%
\begin{eqnarray*}
\mathbb{E}\left[ \mathcal{L}_{\ell }^{2}\right] &\leq &\lim
\inf_{\varepsilon \rightarrow 0}\mathbb{E}\left[ \left\{ \int_{\mathbb{S}%
^{2}}\left\{ \chi _{\varepsilon }(f_{\ell }(x))\left\Vert \nabla f_{\ell
}(x)\right\Vert \right\} dx\right\} ^{2}\right] \\
&=&\lim \inf_{\varepsilon \rightarrow 0}\mathbb{E}\left[ \mathcal{L}_{\ell
;\varepsilon }^{2}\right] \\
&\leq &\lim \sup_{\varepsilon \rightarrow 0}\mathbb{E}\left[ \mathcal{L}%
_{\ell ;\varepsilon }^{2}\right] \\
&=&\lim \sup_{\varepsilon \rightarrow 0}\mathbb{E}\left[ \left\{ \int_{%
\mathbb{S}^{2}}\left\{ \chi _{\varepsilon }(f_{\ell }(x))\left\Vert \nabla
f_{\ell }(x)\right\Vert \right\} dx\right\} ^{2}\right] \\
&=&\lim \sup_{\varepsilon \rightarrow 0}\mathbb{E}\left[ \left\{ \int_{%
\mathbb{R}}\mathcal{L}_{\ell }(u)\chi _{\varepsilon }(u)du\right\} ^{2}%
\right] \\
&\leq &\lim \sup_{\varepsilon \rightarrow 0}\int_{\mathbb{R}}\mathbb{E}\left[
\mathcal{L}_{\ell }^{2}(u)\right] \chi _{\varepsilon }(u)du=\mathbb{E}\left[
\mathcal{L}_{\ell }^{2}\right].
\end{eqnarray*}%
We have thus shown that $\mathbb{E}\left[ \mathcal{L}_{\ell ;\varepsilon
}^{2}\right] \rightarrow \mathbb{E}\left[ \mathcal{L}_{\ell }^{2}\right] ,$
and the proof is complete.

\end{document}